\documentclass{amsproc}

\newtheorem{theorem}{Theorem}[section]
\newtheorem{thm}[theorem]{Theorem}
\newtheorem{lemma}[theorem]{Lemma}
\newtheorem{prop}[theorem]{Proposition}
\newtheorem{remarks}[theorem]{Remarks}

\theoremstyle{definition}

\theoremstyle{remark}
\newtheorem{remark}[theorem]{Remark}

\numberwithin{equation}{section}

\newcommand{\curl}{{\operatorname{curl}\,}}
\newcommand{\Div}{{\operatorname{div}\,}}
\newcommand{\Spec}{\operatorname{Spec}}
\newcommand{\dist}{\operatorname{dist}}

\begin{document}

\title{Optimal Uniform Elliptic Estimates for the Ginzburg-Landau System}

%    Information for first author
\author{S. Fournais}
%    Address of record for the research reported here
\address{CNRS, Laboratoire de Math\'{e}matiques d'Orsay, Orsay Cedex, F-91405;
Univ Paris-Sud, Orsay Cedex, F-91405.}

\email{soeren.fournais@math.u-psud.fr}
\thanks{Both authors were supported by the European Research Network
`Postdoctoral Training Program in Mathematical Analysis of Large
Quantum Systems' with contract number HPRN-CT-2002-00277 and by the ESF
Scientific Programme in Spectral Theory and Partial Differential
Equations (SPECT). Part of this work was carried out while the authors visited the Isaac Newton Institute and its hospitality is gratefully acknowledged.}

%    Information for second author
\author{B. Helffer}
\address{Laboratoire de Math\'{e}matiques d'Orsay, 
Orsay Cedex, F-91405; CNRS, Orsay Cedex, F-91405.
}
\email{bernard.helffer@math.u-psud.fr}
%\thanks{Support information for the second author.}

%    General info
\date{\today}

%\dedicatory{This paper is dedicated to our advisors.}


\begin{abstract}
We reconsider the elliptic estimates for magnetic operators in two and three dimensions used in connection with Ginzburg-Landau theory.
Furthermore we discuss the so-called blow-up technique in order to obtain optimal estimates in the limiting cases.
\end{abstract}

\maketitle

\section{Introduction}

In the analysis of the Ginzburg-Landau system, notably in the study of superconductors of Type II and in the parameter regime known as `above $H_{C_2}$', one often needs to estimate the distance between the induced magnetic vector potential, ${\bf A}$, and the fixed exterior magnetic potential, ${\bf F}$, in various norms. In the literature such estimates are found in varying generality scattered over different publications (c.f. \cite{LuPa1,LuPa2,LuPa3,LuPa5}, \cite{He-Pan},...). 

These estimates come in two types. ~\\
The first set of estimates is deduced from the ellipticity 
of the Ginzburg-Landau system. In this way one obtains the desired estimates in (Sobolev) norms, $W^{s,p}$, for $p< +\infty$
(by imbedding theorems also estimates in H\"{o}lder norms, $C^{s',\alpha}$, $\alpha<1$, are obtained). The challenge here is to get inequalities with the right dependence on the magnetic field strength (as opposed to the vector potential).
This part of the analysis is valid in a large parameter regime and is essentially functional analytical.

The second set of estimates corresponds to the cases $p=\infty$ above and uses the first set of estimates as input. One proves that it is possible to go to these limiting cases essentially without loss in the parameter measuring the magnetic field strength.
These inequalities are asymptotic in the sense that they depend on a certain parameter to be sufficiently large and are valid in a much smaller parameter regime (`above $H_{C_2}$').
The proof of these estimates uses the fact that a natural limiting equation has no non-trivial solutions and the proof is therefore much more intrinsically PDE in spirit. This technique is often called a `blow-up argument' in the literature.

Roughly the first half of this article contains the proof of the basic elliptic regularity results for the Ginzburg-Landau system in $2$ and $3$-dimensions. The main result in $2D$ is stated below as Theorem~\ref{thm:Elliptic}, the corresponding $3D$ result is Theorem~\ref{thm:elliptic-3D}.
The $2D$-result is a slight improvement over the analogous Proposition~3.1 in \cite{LuPa1} and the first motivation for writing this paper was to give a short
rather self-contained  proof of that proposition.
Also the regime of validity of the estimate is clarified, in particular our result is true essentially without condition on the external field. Hence  its domain of validity covers also the region around the second critical field, which can be interesting (see \cite{Panhc2}).

As described above, these elliptic results, i.e.~Theorem~\ref{thm:Elliptic} or \cite[Proposition~3.1]{LuPa1}, are the basic input to control a `blow-up' approach as in \cite[Section~4]{LuPa1} or
\cite[Section~4]{He-Pan}.
In Section~\ref{asymptotic} we describe this approach and give the main results, see Propositions~\ref{prop:1}, \ref{prop:blowup} and ~\ref{prop:PsiLilleInt}.
Here also we obtain slightly more precise versions than previous results on the subject.
Such results have been used in particular cases, for instance in ~\cite{Panhc2} and \cite{AlHe}.

This work was partially motivated by discussions with X-B.~Pan and S.~Serfaty and we thank them for the encouragement.

\noindent{\bf Notation.}\\
We will use the standard Sobolev spaces $W^{s,p}$. Furthermore we will use H\"{o}lder spaces.
Let us fix the definition of the norm in the H\"{o}lder spaces $C^{n,\alpha}$. For a smooth bounded domain $\Omega$, $n \in {\mathbb N}$, $\alpha \in (0,1)$, the space $C^{n,\alpha}(\overline{\Omega})$ is the set of functions $u$ having $n$ H\"{o}lder continuous derivatives in $\overline{\Omega}$ and such that the norm
\begin{align}
\| u \|_{C^{n,\alpha}(\overline{\Omega})} :=
\sum_{|\beta| \leq n } \| \partial^{\beta} u \|_{L^{\infty}(\Omega)} + 
\sum_{|\beta| = n} \sup_{x,y \in \overline{\Omega}} \frac{|\partial^{\beta} u(x) - \partial^{\beta} u(y)|}{|x-y|^{\alpha}},
\end{align}
is finite. In the case where $\alpha=0$, the last sum is omitted.

\section{Integration by parts}
\subsection{The case of dimension $2$}$\,$\\
We will use the following notation for the magnetic derivatives
\begin{align}
{\bf D} = (D_1, D_2) = (-i\nabla + B {\bf A}).
\end{align}
The magnetic Laplacian is now the operator ${\mathcal H} := {\bf D}^2 = D_1^2 + D_2^2$.

\begin{prop}\label{prop:IntParts}~\\
Let $\Omega \subset {\mathbb R}^2$ be a regular bounded domain.
Suppose that $\psi \in W^{2,2}(\Omega)$ satisfies magnetic Neumann boundary conditions
\begin{align}
\label{eq:Neumann}
\nu \cdot D\psi \big |_{\partial \Omega} = 0.
\end{align}
Then
\begin{align}
\sum_{j,k} \| D_j D_k \psi \|_{L^2(\Omega)}^2&= B^2 \int_{\Omega} (\curl {\bf A})^2 |\psi |^2\,dx  + \int_{\Omega} | {\mathcal H} \psi |^2\,dx\nonumber\\
&\quad+
2B \int_{\Omega} (\curl {\bf A}) \Im( D_1 \psi \, \overline{ D_2 \psi })\,dx.
\end{align}
\end{prop}

\begin{remark}~\\
This formula appears in \cite{LuPa2} with an additional boundary term
 that we are able to show to be zero in the case of a (magnetic) Neumann-condition.
\end{remark}

$\,$
\begin{proof}[Proof of Proposition \ref{prop:IntParts}]~\\
The proof consists of a tedious but elementary calculation.
First we calculate without using the boundary condition and on functions in $C^\infty(\overline{\Omega})$.
\begin{align*}
\sum_{j,k} \| &D_j D_k \psi \|_{L^2(\Omega)}^2 =
\sum_{j,k} \Re \int_{\Omega} \overline{ D_j D_k \psi } \; \big( D_k D_j \psi -iB(\partial_j A_k - \partial_k A_j)\psi \big)\,dx \\
&= \sum_{j,k} \Re \{ \int_{\Omega} \overline{ D_k \psi }  \; D_j D_k D_j \psi + 
\int_{\partial \Omega} \overline{ D_k \psi } \; \nu_j \; D_k D_j \psi \,d\sigma \}\\
&\quad - \Re\{ iB \int_{\Omega} (\curl {\bf A}) \psi \; \overline{ D_1 D_2 \psi - D_2 D_1 \psi }\,dx\}
\displaybreak[2]\\
&=
B^2 \int_{\Omega} (\curl {\bf A})^2 |\psi |^2\,dx 
+ \sum_k \Re \int_{\Omega} \overline{ D_k \psi } \; D_k {\mathcal H} \psi \,dx\\
&\quad + \sum_{j,k} \Re \int_{\Omega} \overline{ D_k \psi } \; (-iB) (\partial_j A_k - \partial_k A_j) D_j \psi\,dx\\
&\quad +\sum_{j,k} \Re \int_{\partial \Omega} \overline{ D_k \psi } \; \nu_j \; D_k D_j \psi \,d\sigma  \\
&=
B^2 \int_{\Omega} (\curl {\bf A})^2 |\psi |^2\,dx  + \int_{\Omega} | {\mathcal H} \psi |^2\,dx
+
2B \int_{\Omega} (\curl {\bf A}) \; \Im( D_1 \psi \; \overline{ D_2 \psi })\,dx\\
&\quad +  \Re \int_{\partial \Omega} \big\{ (\nu \cdot \overline{{\bf D} \psi}) \; {\mathcal H} \psi 
+ \sum_{j,k} \overline{ D_k \psi } \; \nu_j \; D_k D_j \psi \big\}\,d\sigma.
\end{align*}
We now apply the Neumann boundary condition. That makes the first boundary term vanish. 
The second boundary term we can rewrite as follows,
\begin{align*}
\Re \int_{\partial \Omega} \sum_{j,k} \overline{ D_k \psi } \; \nu_j \; D_k D_j \psi \,d\sigma
=
\Re(a+b),
\end{align*}
with
\begin{align*}
a&:=  \int_{\partial \Omega} \sum_{j,k} \overline{ D_k \psi } \;  D_k \nu_j  D_j \psi \,d\sigma,\\
b&:=   i \int_{\partial \Omega} \sum_{j,k} \overline{ D_k \psi } \,(\partial_k \nu_j) D_j \psi \,d\sigma\;.
\end{align*}
To analyze $a,b$ we introduce a unit vector $\tau$ parallel to the boundary and define $D_{\tau}:=\tau \cdot {\bf D}$, $D_{\nu}:=\nu \cdot {\bf D}$.

Let us start by proving that $\Re (b)$ vanishes. Taking the real part, we find
\begin{align}
\Re(b)=\frac{i}{2} \int_{\partial \Omega} \langle {\bf D}\psi; M {\bf D}\psi \rangle_{{\mathbb C}^2}\,d\sigma,
\end{align}
where $M$ is the matrix with entries $M_{j,k} = \partial_j \nu_k - \partial_k \nu_j$. 
It clearly suffices to prove that the integrand is real in order to conclude that $\Re (b)=0$.

Writing ${\bf D}\psi = (D_{\tau} \psi )\tau + (D_{\nu} \psi )\nu$ and using the boundary condition, we find that the integrand satisfies
$$
\langle {\bf D}\psi; M {\bf D}\psi \rangle_{{\mathbb C}^2} = | D_{\tau}\psi |^2 \langle \tau ; M \tau \rangle_{{\mathbb C}^2},
$$
which is manifestly real since $M, \tau$ are real.
Thus $\Re (b)=0$.

Using the Neumann boundary condition and the fact that $(\tau, \nu)$ is an orthogonal basis for ${\mathbb C}^2$, we can rewrite $a$ as
\begin{align*}
a= \int_{\partial \Omega} \overline{ D_{\tau} \psi }\; D_{\tau} D_{\nu} \psi \,d\sigma.
\end{align*}
Since (the vector-field part in) $D_{\tau}$ is a derivative along the boundary, and since \break $D_{\nu} \psi \big|_{\partial \Omega} =0$\,, we find
$D_{\tau} D_{\nu} \psi\big|_{\partial \Omega} =0$.
Thus clearly $a$ vanishes.
\end{proof}

We now get an interesting elliptic inequality
for $2D$ magnetic problems with Neumann boundary conditions.

\begin{lemma}
\label{lem:IntParts}~\\
Let $\Omega \subset {\mathbb R}^2$ be a regular domain.
Suppose that $\psi \in C^{\infty}(\overline{\Omega})$ satisfies magnetic Neumann boundary conditions. Then
for all $p_1, p_2 \in [1, +\infty]$ we have
\begin{align}
\label{eq:2D-PartInt}
\sum_{j,k} \| D_j D_k \psi \|_{L^2(\Omega)}^2
&\leq 
3B^2 \| \psi \|_2^2 + 2 \| {\mathcal H} \psi \|_2^2
+2 B^2 \| \curl {\bf A} - 1 \|_{2p_1}^2 \| \psi \|_{2q_1}^2\nonumber\\
&\quad
+ 2B \| \curl {\bf A} - 1 \|_{p_2} \| {\bf D} \psi \|_{2q_2}^2,
\end{align}
where $q_j$ is the conjugate exponent to $p_j$, i.e. $p_j^{-1} + q_j^{-1} = 1$.
\end{lemma}

\begin{proof}~\\
The proof is direct using the identity in Proposition~\ref{prop:IntParts}---replacing $\curl {\bf A}$ by $(\curl {\bf A} - 1) + 1$---and H\"{o}lder's inequality. The term $B \| {\bf D} \psi \|_2^2$ is estimated as
$$
B \| {\bf D} \psi \|_2^2 = B \langle \psi , {\mathcal H} \psi \rangle \leq  B^2 \| \psi \|_2^2 + \| {\mathcal H}\psi \|_2^2, 
$$
where the Neumann boundary condition is used to get the identity.
\end{proof}
\subsection{The case of dimension $3$}~\\
The same calculation as in the $2D$ case yields, using of course that $\psi$ satisfies the
Neumann condition~:
\begin{align}
\label{eq:PartInt3D}
\sum_{j,k} \| D_j D_k \psi \|_{L^2(\Omega)}^2&= B^2 \int_{\Omega}
(\curl {\bf A})^2 |\psi |^2\,dx  +   \int_{\Omega} | {\mathcal H} \psi |^2\,dx\nonumber\\
&\quad+
2B\int_{\Omega} (\curl {\bf A}) \Im\left(
\begin{matrix}
D_2\psi \, \overline{D_3 \psi} \\
D_3\psi \, \overline{D_1 \psi} \\
D_1\psi \, \overline{D_2 \psi}
\end{matrix}
\right)\,dx
 +\Re  b\,,
\end{align}
with 
\begin{align}
b:=  i \int_{\partial \Omega} \sum_{j,k} \overline{ D_k \psi } \,(\partial_k \nu_j) \,D_j \psi \,d\sigma\;,\quad\quad
\curl {\bf A} = \left(\begin{matrix} \partial_2 A_3 - \partial_3 A_2 \\
\partial_3 A_1 - \partial_1 A_3 \\ \partial_1 A_2 - \partial_2 A_1
\end{matrix}\right)
\end{align}
In the $3$ dimensional case we are not able to prove that $b$ vanishes,
but this boundary term can be controlled as follows by trace theorems.

Since the derivatives of $\nu$ are bounded we can estimate
$$
| b | \leq C \| {\bf D} \psi \|_{L^2(\partial \Omega)}^2. 
$$
Notice that elementary identity
$$
|u(0)|^2 = -2 \int_0^{\infty} \frac{d}{dt} |u(t)|^2\,dt ,
$$
for $u\in H^1(\overline{{\mathbb R}_+})$, implies the inequality
$$
|u(0)| \leq \sqrt{2} \, \| u \|_{L^2({\mathbb R}_+)} \| u' \|_{L^2({\mathbb R}_+)}.
$$
Implementing this inequality in a suitable set of coordinates near the boundary, one sees that there exists a constant $C>0$, such that for all $\epsilon <1$ and all $f \in W^{1,2}(\Omega)$ we have
\begin{align*}
\| f \|_{L^2(\partial \Omega)}^2 & \leq C \epsilon^{-1} \| f \|_{L^2(\Omega)}^2 + \epsilon \| f \|_{W^{1,2}(\Omega)}^2.
\end{align*}

We will choose $\epsilon = 1/2$ and apply the resulting inequality to $f = |{\bf D} \psi|$. Combining with \eqref{eq:PartInt3D} we thereby get
\begin{align}
\label{eq:IntParts3D}
\sum_{j,k} \| D_j D_k \psi \|_{L^2(\Omega)}^2
&\leq
C \| {\bf D} \psi \|_{L^2(\Omega)}^2 
+
B^2 \int_{\Omega}
(\curl {\bf A})^2 |\psi |^2\,dx  +   \int_{\Omega} | {\mathcal H} \psi |^2\,dx\nonumber\\
&\quad+
2B\Big| \int_{\Omega} (\curl {\bf A}) \Im\left(
\begin{matrix}
D_2\psi \, \overline{D_3 \psi} \\
D_3\psi \, \overline{D_1 \psi} \\
D_1\psi \, \overline{D_2 \psi}
\end{matrix}
\right)\,dx\Big|.
\end{align}
Let $\beta=(0,0,1)$ denote the unit constant magnetic field. 
The $3D$ result analogous to Lemma~\ref{lem:IntParts} is the following.

\begin{lemma}
\label{lem:Intparts-3D}~\\
Let $\Omega \subset {\mathbb R}^3$ be a smooth domain with compact boundary. 
Then there exists a constant $C>0$ such that for all $p_1,p_2 \in [1,+\infty]$ and all 
$\psi \in C^{\infty}(\overline{\Omega})$ satisfying magnetic Neumann conditions, we have
\begin{align}
\sum_{j,k} \| D_j &D_k \psi \|_{L^2(\Omega)}^2\nonumber\\
&\leq 
C \Big\{B^2 \| \psi \|_2^2 + 
(1+B) \| {\bf D} \psi \|_2^2
+ \| {\mathcal H} \psi \|_2^2
\nonumber\\
&\quad\quad\quad
+B^2 \| \curl {\bf A} - \beta \|_{2p_1}^2 \| \psi \|_{2q_1}^2
+ B \| \curl {\bf A} - \beta \|_{p_2} \| {\bf D} \psi \|_{2q_2}^2\Big\},
\end{align}
where $q_j$ is the conjugate exponent to $p_j$, i.e. $p_j^{-1} + q_j^{-1} = 1$.
\end{lemma}

\section{Regularity for the solutions  of the  Ginzburg-Landau system }
\subsection{The $2D$ case}$\,$\\
We recall that the Ginzburg-Landau functional is given by
\begin{multline}
\label{eq:GL_F}
{\mathcal E}[\psi,{\bf A}] = {\mathcal
E}_{\kappa,H}[\psi,{\bf A}]  =
\int_{\Omega} \Big\{ |p_{\kappa H {\bf A}}\psi|^2 
- \kappa^2|\psi|^2
+\frac{\kappa^2}{2}|\psi|^4 \\
+ \kappa^2 H^2
|\curl {\bf A} - 1|^2\Big\}\,dx\;,
\end{multline}
with
$ (\psi, {\bf A})  \in W^{1,2}(\Omega;{\mathbb C})\times
W^{1,2}(\Omega;{\mathbb R}^2)$. 
We have introduced the notation $p_{\bf A}$ for the operator $(-i\nabla + {\bf A})$.

Let us fix the choice of (London) gauge by imposing that 
\begin{align}
\label{eq:gauge}
\Div {\bf A} &= 0 \quad \text{ in } \Omega\;, & {\bf A} \cdot \nu = 0 \quad \text{ on } \partial \Omega\;.
\end{align}

We also recall that a minimizer of the Ginzburg-Landau functional
 satisfies the  Ginzburg-Landau equations.
\begin{subequations}
\label{eq:GL}
\begin{align}
\left.\begin{array}{c}
p_{\kappa H {\bf A}}^2\psi =
\kappa^2(1-|\psi|^2)\psi \\
\label{eq:equationA}
\curl^2 {\bf A} =-\tfrac{i}{2\kappa H}(\overline{\psi} \nabla
\psi - \psi \nabla \overline{\psi}) + |\psi|^2 {\bf A}
\end{array}\right\} &\quad \text{ in } \quad \Omega \, ;\\
\left. \begin{array}{c}
(p_{\kappa H {\bf A}} \psi) \cdot \nu = 0 \\
\curl {\bf A} - 1 = 0
\end{array} \right\} &\quad \text{ on } \quad \partial\Omega \, .
\end{align}
\end{subequations}
Here $\curl (A_1,A_2) = \partial_{x_1}A_2 -
\partial_{x_2}A_1$, and
$$
\curl^2 {\bf A} =
(\partial_{x_2}(\curl {\bf A}),-\partial_{x_1}(\curl {\bf
A})) \, .
$$

When discussing the Ginzburg-Landau system \eqref{eq:GL} we will always impose the gauge condition \eqref{eq:gauge}.

For a solution of the Ginzburg-Landau system, we deduce the following 
standard estimate \cite{DGP,Giorgi-Phillips} ~:
\begin{equation}\label{infini}
\|\psi\|_\infty \leq 1\;.
\end{equation}

Furthermore, we let ${\bf F}$ denote the vector potential generating a constant magnetic field, which more explicitly satisfies~:
\begin{align}
\curl {\bf F} = 1 \quad \text{ in } \Omega,\quad 
\Div {\bf F} = 0 \quad \text{ in } \Omega,\quad 
{\bf F} \cdot \nu\big|_{\partial \Omega} = 0.
\end{align}

\begin{thm}\label{thm:Elliptic}~\\
Let $\Omega \subset {\mathbb R}^2$ be a smooth, bounded domain.
There exist a constant $C$, and, for any $\alpha \in(0,1)$ and $p\in (1,+\infty)$, 
constants $\widehat C_\alpha$ and $\widetilde C_p$,  such that,
if $(\psi,{\bf A})$ is any solution of the Ginzburg-Landau system \eqref{eq:GL} with parameters $\kappa, H>0$, then
\begin{align}
\label{eq:DD1}
&\sum_{j,k} \| D_j D_k \psi \|_{L^2(\Omega)} \leq C (1+ \kappa H +\kappa^2) \| \psi \|_2,\\
\label{eq:CAF1}
& \|\curl  {\bf A} - 1 \|_{C^{0,\alpha}(\overline{\Omega})} 
\leq \widehat C_\alpha  \frac{1+ \kappa H +\kappa^2}{\kappa H} \| \psi \|_2 \| \psi \|_{\infty},\\
and & \nonumber \\
& \| \curl {\bf A} - 1 \|_{W^{1,p}(\Omega)}\leq \widetilde C_p  \frac{1+ \kappa H +\kappa^2}{\kappa H} \| \psi \|_2 \| \psi \|_{\infty}\;.\label{eq:CAF2}
\end{align}
\end{thm}

\begin{remarks}~
\begin{itemize}
\item Using the $W^{k,p}$-regularity of the Curl-Div system (see \cite{AgDoNi2}, see also \cite{Tem} for the case $p=2$),
we obtain from \eqref{eq:CAF2} the estimate
\begin{equation}\label{eq:A2}
 \| {\bf A} - {\bf F} \|_{W^{2,p}(\Omega)}\leq \widetilde D_p  \frac{1+ \kappa H +\kappa^2}{\kappa H}  \| \psi \|_2 \| \psi \|_{\infty}\;.
\end{equation}
Hence, using the Sobolev injection Theorem,
\begin{equation}\label{eq:AF1}
\| {\bf A} - {\bf F} \|_{C^{1,\alpha}(\overline{\Omega})} 
\leq \widehat D_\alpha 
\frac{1+ \kappa H +\kappa^2}{\kappa H}  \| \psi \|_2 \| \psi \|_{\infty},
\end{equation}
for all $\alpha \in [0,1)$.
\item In the applications, $H$ is of the same order as $\kappa$, so
 \eqref{eq:AF1}
 gives that $({\bf A}-{\bf F})$ is uniformly bounded in $C^{1,\alpha}(\overline{\Omega})$
 in this regime, for any $\alpha <1$.
\item
We have in particular obtained a complete proof of the basic Proposition~3.1
 in \cite{LuPa1} with actually an improvement of the right hand side
 and an extension of the regime of parameters $(\kappa,H)$ for which the
 estimate is true.
\item When in addition, $\frac{\kappa}{H} \geq 1 +b$ (with $b>0$), 
 V. Bonnaillie-No\"el and S. Fournais have given in \cite{BonFo} a very simple proof (in comparison with \cite{He-Pan} or \cite{FournaisHelffer3}) showing
 that for a minimizer $(\psi,{\bf A})$ of the Ginzburg-Landau functional, one has for some constants $C_b,\kappa_b>0$,
\begin{align}
\label{eq:Decay}
\|\psi\|_2 \leq C \kappa^{-\frac 12}\; \|\psi\|_\infty\;,
\end{align}
for all $\kappa\geq \kappa_b$.
It is important to note that the proof in \cite{BonFo} does not use the elliptic estimates discussed in the present paper. This is in contradistinction to previous derivations of inequalities like \eqref{eq:Decay} (see \cite{He-Pan}) which use \eqref{eq:AF1} as an input.
\item With almost no modification we can treat the case when $\mathcal H_0= \curl {\bf} F$
 is a regular function in $\overline{\Omega}$, instead of a constant. The second equation in
 \eqref{eq:GL} becomes in this case
 $$
\curl\left(\curl  {\bf A} - \mathcal H_0\right)  =-\tfrac{i}{2\kappa H}(\overline{\psi} \nabla
\psi - \psi \nabla \overline{\psi}) + |\psi|^2 {\bf A}  \quad \text{ in } \quad \Omega \, .
$$
\end{itemize}
\end{remarks}

$\,$

\begin{proof}[Proof of Theorem~\ref{thm:Elliptic}]~\\
Recall the estimate \eqref{infini}.
Furthermore,
multiplying \eqref{eq:equationA} by $\overline{\psi}$, integrating and implementing \eqref{infini}, we obtain
\begin{equation}\label{cine}
\| (-i\nabla + \kappa H {\bf A})\psi \|_2 \leq \kappa \| \psi \|_2\;.
\end{equation}

Using the second equation of the G-L system, we get 
\begin{equation}\label{ellp1}
\| \nabla (\curl {\bf A}-1) \|_p \leq  \frac{C}{\kappa H}
 \|\psi\|_\infty\; \| (-i\nabla + \kappa H {\bf A})\psi \|_p\;.
\end{equation}
But using the property that $\curl {\bf A} -1$ satisfies the Dirichlet
 condition, this implies
\begin{equation}\label{ellp2}
\|  \curl {\bf A}-1 \|_p\leq  \frac{C}{\kappa H}
 \|\psi\|_\infty\; \| (-i\nabla + \kappa H {\bf A})\psi \|_p\;.
\end{equation}
When $p=2$, we can then implement the control of $\| (-i\nabla + \kappa H {\bf A})\psi \|_2$ obtained in \eqref{cine} and get
\begin{equation}\label{ineqimproved}
\|  \curl {\bf A}-1 \|_2\leq  \frac{C}{ H}
 \|\psi\|_\infty\; \|\psi \|_2\;.
\end{equation}
Note also that we actually get the stronger estimate 
\begin{equation*}
\|  \curl {\bf A}-1 \|_{H^1} \leq  \frac{C}{ H}
 \|\psi\|_\infty\; \|\psi \|_2\;,
\end{equation*}
but this will not be used.

Using then \eqref{ellp1}, \eqref{ellp2} and the Sobolev-injection Theorem,
we get that for any $\alpha \in (0,1)$ there exists $p=p(\alpha)$ such that
\begin{align}\label{eq:C-1-alphab}
\|\curl {\bf A} - 1 \|_{C^{0,\alpha}(\overline{\Omega})}
&\leq   \frac{C''}{\kappa H} \| \psi \|_{\infty} \| (-i\nabla + \kappa H {\bf A})\psi \|_{L^{p}(\Omega)}\;. \end{align}
Here all the constants depend on $\alpha$. 
We now use the (pointwise) diamagnetic inequality
\begin{equation}
|\, \nabla |\chi|\,|\leq |(\nabla + i\kappa A) \chi|\;,
\end{equation}
(actually applied with $\chi= (-i\partial_k +\kappa H A_k)\psi$)
in order to get
\begin{align}\label{inegaga}
\|\curl { \bf A} - 1\|_{C^{0,\alpha}(\overline{\Omega})}
&\leq  \frac{C''' }{\kappa H} \| \psi \|_{\infty} \sum_{j,k} \| (-i\partial_j + \kappa H A_j)(-i\partial_k + \kappa H A_k)\psi\|_2 \nonumber\\
& \quad +   \frac{C'''}{\kappa H}\| \psi \|_{\infty} \sum_k 
\|(-i\partial_k + \kappa H A_k)\psi\|_2 \;.
\end{align}

We will insert the estimates obtained above in \eqref{eq:2D-PartInt}.
In preparation for this, we estimate using \eqref{cine},
\begin{align}
\label{eq:Calpha-curl}
\| \curl {\bf A} -1 \|_{\infty} \| {\bf D} \psi\|_2^2 
\leq \kappa^2 \|\curl  {\bf A} - 1 \|_{C^{0,\alpha}(\overline{\Omega})} \| \psi \|_2^2.
\end{align}
Also, using \eqref{ineqimproved}, we obtain
\begin{align}
\label{eq:L2-curl}
\| \psi \|_{\infty}^2 \| \curl {\bf A} -1 \|_2^2 
\leq \frac{1}{H^2} \,\| \psi \|_{\infty}^2\, \| \psi \|_2^2.
\end{align}
Using Lemma~\ref{lem:IntParts} with $B = \kappa H$, $p_1=1$ and $p_2=+\infty$, combined with \eqref{eq:equationA}, \eqref{infini}, \eqref{eq:Calpha-curl} and \eqref{eq:L2-curl} we get
\begin{align}
\sum_{j,k} \| D_j D_k \psi \|_{L^2(\Omega)}^2&\leq
3  \{ 1 + (\kappa H)^2 + \kappa^2 + \kappa^4\} \| \psi \|_2^2 
\nonumber\\
&+
 2 \kappa^3 H \|\curl {\bf A} - 1\|_{C^{0,\alpha}}  \| \psi \|_2^2.
\end{align}
Thus,
\begin{align}
\sum_{j,k} \| D_j D_k \psi \|_{L^2(\Omega)}&\leq
C  \Big\{ (1 + \kappa H + \kappa^2) \| \psi \|_2 +
\kappa^{3/2} H^{1/2}\|\curl  {\bf A} - 1 \|_{C^{0,\alpha}}^{1/2} \| \psi \|_2\Big\}\,.\nonumber
\end{align}
Hence, there exists a constant $C'$, such that, for all $\epsilon >0$, we have
\begin{align}\label{eq:DD}
\sum_{j,k} \| D_j D_k \psi \|_{L^2(\Omega)}
&\leq
C'  \Big\{ \frac{1 + \kappa H + \kappa^2}{\kappa H}  \| \psi \|_2 + \epsilon^{-1} \kappa^2 \| \psi\|_2^2 \Big\}
\nonumber\\
&\quad\quad+
\epsilon (\kappa H) \| \curl {\bf A} - 1 \|_{C^{0,\alpha}}\;.
\end{align}
We insert \eqref{eq:DD} and \eqref{cine} in \eqref{inegaga} and  find, for some constant $C>0$,
$$
(1- C \epsilon  \|\psi\|_\infty )\,\|\curl  {\bf A} - 1\|_{C^{0,\alpha}}
 \leq C \|\psi_\infty\| \frac{1 + \kappa H + \kappa^2}{\kappa H} \|\psi\|_2
 + C \epsilon^{-1} \frac{\kappa}{H} \|\psi\|^2_2\;.
$$
Taking $\epsilon =\frac{1}{2C}$ and using \eqref{ineqimproved} leads to the expected
 \begin{align}
\label{eq:AF}
\|\curl {\bf A} - 1 \|_{C^{0,\alpha}} 
\leq \widehat C\, \frac{1 + \kappa H + \kappa^2}{\kappa H}  \| \psi \|_2 \| \psi \|_{\infty},
\end{align}
where the constant $\widehat C$ is independent of $\kappa$ and $H$.

Inserting the bound \eqref{eq:AF} in \eqref{eq:DD} yields
\begin{align}
\label{eq:DD2}
\sum_{j,k} \| D_j D_k \psi \|_{2} \leq C (1 + \kappa H + \kappa^2) \| \psi \|_2.
\end{align}
The proof of the last statement in Theorem~\ref{thm:Elliptic} is obtained by starting again
from the right inequality in \eqref{eq:C-1-alphab}, and implementing
\eqref{eq:CAF1}.
This finishes the proof of Theorem~\ref{thm:Elliptic}.
\end{proof}

\subsection{The $3D$ case}$\,$\\
In three dimensions a generalization of the Ginzburg-Landau functional 
which is often considered
is 
\begin{multline}
\label{eq:GL_F-3D}
{\mathcal E}_{\kappa,H}^{3D}[\psi,{\bf A}]  =
\int_{\Omega} \Big\{ |p_{\kappa H {\bf A}}\psi|^2 
- \kappa^2|\psi|^2
+\frac{\kappa^2}{2}|\psi|^4 \Big\}\,dx\\
+ \kappa^2 H^2 \int_{{\mathbb R}^3}
|\curl {\bf A} - \beta|^2 \,dx\;,
\end{multline}
where $\beta$ is the external magnetic field. We will choose $\beta = (0,0,1)$ i.e. a constant magnetic field, but more general situations could easily be considered.
Notice that the field integral is over ${\mathbb R}^3$ instead of $\Omega$.

We will consider the case of smooth, bounded $\Omega$.
Let $\dot{H}^1({\mathbb R}^3)$ denote the homogeneous Sobolev space, i.e. the closure of $C_0^{\infty}({\mathbb R}^3)$ under the norm $u \mapsto \| u \|_{\dot{H}^1({\mathbb R^3})} :=\| \nabla u \|_{L^2({\mathbb R}^3)}$.
Let furthermore ${\bf F}$ denote the vector potential generating the constant magnetic field, 
${\bf F}(x) = (-x_2/2, x_1/2, 0)$. 
Clearly, $\Div {\bf F}=0$.
Then the natural variational space for the functional ${\mathcal E}_{\kappa,H}^{3D}$ is $(\psi, {\bf A}) \in W^{1,2}(\Omega) \times \dot{H}^1_{{\rm div}, {\bf F}}$, where
$$
\dot{H}^1_{{\rm div}, {\bf F}} = \{ {\bf A} \;:\; \Div {\bf A} = 0, \text{ and } {\bf A}- {\bf F} \in \dot{H}^1({\mathbb R}^3) \}.
$$

Minimizers of ${\mathcal E}_{\kappa,H}^{3D}$ are weak solutions of the Euler-Lagrange equations
\begin{subequations}
\label{eq:GL-3D}
\begin{align}
\label{eq:equationpsi-3D}
&p_{\kappa H {\bf A}}^2\psi = \kappa^2(1-|\psi|^2)\psi \quad &&\text{ in } \quad \Omega\\
\label{eq:equationA-3D}
&\curl^2 {\bf A} =-\tfrac{1}{\kappa H} \Re(\overline{\psi} p_{\kappa H {\bf A}}\psi) 1_{\Omega}
\quad &&\text{ in } \quad {\mathbb R}^3\\
&(p_{\kappa H {\bf A}} \psi) \cdot \nu = 0, &&\text{ on } \quad \partial\Omega \, .
\end{align}
\end{subequations}
As in the $2D$ case we will give estimates valid for general solutions $(\psi, {\bf A}) \in W^{1,2}(\Omega) \times \dot{H}^1_{{\rm div}, {\bf F}}$ of \eqref{eq:GL-3D} not only for minimizers of ${\mathcal E}_{\kappa,H}^{3D}$.
As we will see below, the fact that we do not have a boundary condition for ${\bf A}$ will both be a simplification and a complication.
The following $3D$ result is similar to \cite[Lemma~3.3]{PanJMP}.

\begin{thm}\label{thm:elliptic-3D}~\\
Let $\Omega\subset {\mathbb R}^3$ be a smooth, bounded domain.
For all $\alpha <1/2$ and all $1\leq p\leq 6$ there exist constants $C_{\alpha}, C_p$ such that
for all $\kappa, H >0$, and all solutions 
$(\psi, {\bf A}) \in W^{1,2}(\Omega) \times \dot{H}^1_{{\rm div}, {\bf F}}$ of \eqref{eq:GL-3D},
\begin{align}
\| {\bf A}- {\bf F} \|_{W^{2,p}(\Omega)}
&\leq
C_p \frac{1 + \kappa H + \kappa^2}{\kappa H} \| \psi \|_{\infty}\| \psi \|_{L^2(\Omega)},\\
\| {\bf A}- {\bf F} \|_{C^{1,\alpha}(\overline{\Omega})}
&\leq
C_{\alpha} \frac{1 + \kappa H + \kappa^2}{\kappa H}\| \psi \|_{\infty}\| \psi \|_{L^2(\Omega)}.
\end{align}
\end{thm}

\begin{proof}~\\
The equations \eqref{infini} and \eqref{cine} remain true for solutions to \eqref{eq:GL-3D}, i.e.
\begin{equation}\label{infini_3D}
\|\psi\|_\infty \leq 1\;,
\end{equation}
and
\begin{equation}\label{cine-3D}
\| (-i\nabla + \kappa H {\bf A})\psi \|_2 \leq \kappa \| \psi \|_2\;.
\end{equation}

We start by noticing that for vector fields $\alpha$ in three dimensions 
the norm $\| \alpha \|_{\dot{H}^1({\mathbb R}^3)}$ is equivalent to the norm
$\| \curl \alpha \|_{L^2({\mathbb R}^3)} + \| \Div \alpha \|_{L^2({\mathbb R}^3)}$.
Furthermore, by the homogeneous Sobolev inequality, the $\dot{H}^1$ norm controls the $L^6$ norm, i.e.~there exists a constant $S_3$ such that
$$
\| u \|_{L^6({\mathbb R}^3)} \leq S_3 \| \nabla u \|_{L^2({\mathbb R}^3)},\quad\quad
\forall u \in C_0^{\infty}({\mathbb R}^3).
$$
Combining these two facts with $\Div( {\bf A} - {\bf F} )= 0$, we find that
\begin{align}
\label{eq:L6-curl}
\| {\bf A} - {\bf F} \|_{L^6({\mathbb R}^3)} \leq C \| \curl {\bf A} - \beta \|_{L^2({\mathbb R}^3)}.
\end{align}

Since $\Div ({\bf A} - {\bf F}) = 0$, the equation \eqref{eq:equationA-3D} can be reformulated as
\begin{align}
\label{eq:equationA-3D-Delta}
\Delta ({\bf A}- {\bf F}) =-\tfrac{1}{\kappa H} \Re(\overline{\psi} \, p_{\kappa H {\bf A}}\psi) 1_{\Omega}
\quad &&\text{ in } \quad {\mathbb R}^3.
\end{align}
Let $B(0,R)$ be the open ball of radius $R$ around the origin.
Elliptic regularity for the Laplacian (see \cite[Theorem~9.11]{GilbargTrudinger}) thus implies for all $p' \in [1, \infty)$, $R>0$, the existence of a constant $C_{p'}(R)$ such that
\begin{align*}
\| {\bf A}- {\bf F} \|_{W^{2,p'}(B(0,R))} \leq C_{p'}(R) \big( \| {\bf A}- {\bf F} \|_{L^{p'}(B(0,2R))} +
\frac{1}{\kappa H} \| \psi \|_{\infty} \| p_{\kappa H {\bf A}}\psi\|_{L^{p'}(\Omega)}\big).
\end{align*}
In particular, for $p'\leq 6$, we can apply the estimate  \eqref{eq:L6-curl} and the compactness of $\overline{B(0,2R)}$ to get
\begin{align}
\| {\bf A}- &{\bf F} \|_{W^{2,p'}(B(0,R))} \nonumber\\
&\leq C_{p'}' (R)\big(  \| \curl {\bf A} - \beta \|_{L^2({\mathbb R}^3)} +
\frac{1}{\kappa H} \| \psi \|_{\infty} \| p_{\kappa H {\bf A}}\psi\|_{L^{p'}(\Omega)}\big),
\end{align}
for all $p'\leq 6$.

Let $R$ be chosen so big that $\Omega \subset B(0,R-1)$. Using once again elliptic regularity and the Sobolev imbedding theorem we find for any $p \in [1,\infty)$,
\begin{align}
\label{eq:Rest-p}
\| {\bf A}- {\bf F} \|_{W^{2,p}(\Omega)} &\leq
C \big( \| {\bf A}- {\bf F} \|_{L^p(B(0,R))} + \frac{1}{\kappa H} \| \psi \|_{\infty} \| p_{\kappa H {\bf A}}\psi\|_{L^{p}(\Omega)}\big)\nonumber\\
&\leq 
C \big( \| {\bf A}- {\bf F} \|_{W^{2,2}(B(0,R))} + \frac{1}{\kappa H} \| \psi \|_{\infty} \| p_{\kappa H {\bf A}}\psi\|_{L^{p}(\Omega)}\big)\nonumber\\
&\leq 
C \big(\| \curl {\bf A} - \beta \|_{L^2({\mathbb R}^3)} +
\frac{1}{\kappa H} \| \psi \|_{\infty} \| p_{\kappa H {\bf A}}\psi\|_{L^{2}(\Omega)} \nonumber\\
&\quad\quad\quad+ \frac{1}{\kappa H} \| \psi \|_{\infty} \| p_{\kappa H {\bf A}}\psi\|_{L^{p}(\Omega)}\big)
\end{align}

Multiplying \eqref{eq:equationA-3D} by ${\bf A}- {\bf F}$ and integrating by parts yields
\begin{align*}
\| \curl {\bf A} - \beta \|^2_{L^2({\mathbb R}^3)}
&=
-\int_{\Omega} ({\bf A}- {\bf F}) \tfrac{1}{\kappa H} \Im(\overline{\psi} p_{\kappa H {\bf A}}\psi)\,dx \\
&\leq 
\frac{1}{\kappa H}
\| {\bf A}- {\bf F} \|_{L^2(\Omega)} \| \psi \|_{\infty} \| p_{\kappa H {\bf A}}\psi\|_{L^2(\Omega)} \\
&\leq
 \frac{C}{\kappa H}
\| {\bf A}- {\bf F} \|_{L^6(\Omega)} \| \psi \|_{\infty} \| p_{\kappa H {\bf A}}\psi\|_{L^2(\Omega)}.
\end{align*}
Implementing the estimates \eqref{eq:L6-curl} and \eqref{cine-3D} we obtain
\begin{align}
\label{eq:curl-beta}
\| \curl {\bf A} - \beta \|_{L^2({\mathbb R}^3)}\leq C H^{-1} \| \psi \|_{\infty} \| \psi \|_{L^2(\Omega)}.
\end{align}
Thus \eqref{eq:Rest-p} becomes, using again \eqref{cine-3D},
\begin{align}
\| {\bf A}- {\bf F} \|_{W^{2,p}(\Omega)} &\leq
C\big( H^{-1} \| \psi \|_{\infty} \| \psi \|_{L^2(\Omega)} + 
\frac{1}{\kappa H} \| \psi \|_{\infty} \| p_{\kappa H {\bf A}}\psi\|_{L^{p}(\Omega)}\big).
\end{align}
By Sobolev imbeddings, the diamagnetic inequality and \eqref{cine-3D}, we therefore find for $p\leq 6$ the estimate
\begin{align}
\label{eq:ContByDD}
\| {\bf A}- {\bf F} \|_{W^{2,p}(\Omega)} &\leq
C\big( H^{-1} \| \psi \|_{\infty} \| \psi \|_{L^2(\Omega)} + 
\frac{1}{\kappa H} \| \psi \|_{\infty} \big\| |p_{\kappa H {\bf A}}\psi|\big\|_{W^{1,2}(\Omega)}\big) \nonumber\\
&\leq
C' \big( H^{-1} \| \psi \|_{\infty} \| \psi \|_{L^2(\Omega)} + 
\frac{1}{\kappa H} \| \psi \|_{\infty} \sum_{j,k} \| D_j D_k \psi\|_{L^2(\Omega)}\big) .
\end{align}
We will use \eqref{eq:ContByDD} for $p=6$ and the Sobolev inequality
\begin{align}
\label{eq:SobIgen}
\| \curl {\bf A} - \beta \|_{L^{\infty}(\Omega)}\leq C \| {\bf A}- {\bf F} \|_{W^{2,6}(\Omega)}.
\end{align}
We use Lemma~\ref{lem:Intparts-3D} with $p_1=1$ and $p_2=+\infty$, and find, by implementing \eqref{eq:equationpsi-3D}, \eqref{infini_3D}, \eqref{cine-3D} and \eqref{eq:curl-beta},
\begin{align*}
\sum_{j,k} \| D_j &D_k \psi\|_{L^2(\Omega)}^2 \\
&\leq
C'\Big\{
(1 +  (\kappa H)^2 + \kappa^4 ) \| \psi \|_{L^2(\Omega)}^2
+ \kappa^3 H \| \curl {\bf A} - \beta \|_{L^{\infty}(\Omega)}\| \psi \|_{L^2(\Omega)}^2
\Big\}.
\end{align*}
So (using $\| \psi \|_{\infty} \leq 1$), for all $\epsilon>0$, 
\begin{align}
\label{eq:FoerPP}
\sum_{j,k} \| D_j D_k \psi\|_{L^2(\Omega)}
\leq
C \Big\{
(1 +  (\kappa H) + &\kappa^2 + \epsilon^{-1} \kappa^2) \| \psi \|_{L^2(\Omega)}\nonumber\\
&+
\epsilon (\kappa H) \| \curl {\bf A} - \beta \|_{L^{\infty}(\Omega)} \Big\}.
\end{align}
Choosing $\epsilon$ sufficiently small and inserting \eqref{eq:SobIgen} and \eqref{eq:FoerPP}
in \eqref{eq:ContByDD} with $p=6$,
we get
\begin{align}
\| {\bf A}- {\bf F} \|_{W^{2,6}(\Omega)}
\leq
C \frac{1 + \kappa^2 + \kappa H}{\kappa H} \| \psi \|_{\infty}\| \psi \|_{L^2(\Omega)}.
\end{align}
Using now Sobolev imbeddings we have proved Theorem~\ref{thm:elliptic-3D}.
\end{proof}

\section{Asymptotic estimates}
\label{asymptotic}
\subsection{Nonexistence of solutions to certain partial differential equations}
\label{nonexist}~\\
We will use the notation ${\bf \tilde F}$ for any vector potential on ${\mathbb R}^2$ or on the half-space
${\mathbb R}^2_+ :=\{(x_1,x_2) \in {\mathbb R}^2 \,\big|\, x_1 > 0\}$ satisfying
$\curl {\bf \tilde F} =1$.

The natural self-adjoint extension of the differential operator $(-i\nabla +{\bf \tilde F})^2$ on $L^2({\mathbb R}^2)$ is known to have spectrum, 
$$
\Spec (-i\nabla +{\bf \tilde F})^2_{L^2({\mathbb R}^2)} = \big\{ 2j + 1, j \in {\mathbb N}\cup \{ 0 \}\big\}.
$$
We also consider the Neumann-realization ${\mathcal H}$ of the same operator but restricted to the half-space
${\mathbb R}^2_+$. This is the operator $(-i\nabla +{\bf \tilde F})^2$ with domain 
$$
\big\{ \psi \in L^2({\mathbb R}^2_+) \,\big|\, (-i\nabla +{\bf \tilde F})^2\psi \in L^2({\mathbb R}^2_+) \quad\text{ and }\quad
\nu \cdot (-i\nabla +{\bf \tilde F}) \psi \big|_{\partial {\mathbb R}_{+}^2 } = 0\big\}.
$$
We define a real number $\Theta_0$ by
\begin{align}\label{eq:InfSpec}
\Theta_0:=\inf \Spec {\mathcal H}.
\end{align}
The number $\Theta_0$ (also denoted by $\beta_0$ by some authors) plays an important role in the analysis of the Ginzburg-Landau system (see for instance \cite{LuPa1,PiFeSt,He-Mo} for information on this spectral constant). Here we will only recall the basic property that  $\Theta_0\in(0,1)$. 

In this subsection we will consider the following PDEs.
\begin{align}
\label{eq:Schr2}
&(-i\nabla +{\bf \tilde F})^2 \psi = \lambda \psi, &&\text{ on } {\mathbb R}^2, \text{ with } \lambda < 1, \\
\label{eq:NL2}
&(-i\nabla +{\bf \tilde F})^2 \psi = \lambda (1- S^2|\psi|^2) \psi, &&\text{ on } {\mathbb R}^2, \text{ with } 0\leq\lambda \leq 1,\\
\label{eq:Schr2+}
&(-i\nabla +{\bf \tilde F})^2 \psi = \lambda \psi, &&\text{ on } {\mathbb R}_{+}^2, \text{ with } \lambda < \Theta_0, \\
\label{eq:NL2+}
&(-i\nabla +{\bf \tilde F})^2 \psi = \lambda (1- S^2|\psi|^2) \psi, &&\text{ on } {\mathbb R}_{+}^2, \text{ with } 0\leq\lambda \leq \Theta_0.
\end{align}
The last two equations \eqref{eq:Schr2+}, \eqref{eq:NL2+} are considered with Neumann boundary condition, i.e.
$\nu \cdot (-i\nabla +{\bf \tilde F}) \psi \big|_{\partial {\mathbb R}_{+}^2 } = 0$.
In order for this boundary condition to be well-defined we assume that $\psi \in H^2_{\rm loc}({\mathbb R}_{+}^2)$.
Also, we assume that the parameter $S\geq 0$ in \eqref{eq:NL2} verifies $S \neq 0$ when $\lambda = 1$, and
similarly, the parameter $S\geq 0$ in \eqref{eq:NL2+}  satisfies $S \neq 0$ when $\lambda = \Theta_0$.

\begin{prop}\label{prop:NonExist}~\\
Let $(\psi, \lambda)$ be a solution to one of the equations \eqref{eq:Schr2}, \eqref{eq:NL2}, \eqref{eq:Schr2+} or \eqref{eq:NL2+} with $\lambda$ in the indicated interval and $\psi$ being globally bounded. Then $\psi =0$.
\end{prop}

\begin{proof}
We only consider the cases on ${\mathbb R}^2_{+}$ since the other statements follow by the same arguments.

Let ${\mathcal H}$ be the operator $(-i\nabla +{\bf \tilde F})^2$ with the Neumann boundary condition. 
We will prove that a non-zero bounded solution to \eqref{eq:Schr2+} or \eqref{eq:NL2+} will provide a contradiction to \eqref{eq:InfSpec} through the variational principle.

Let $\psi \in L^{\infty}({\mathbb R}^2_{+})\setminus \{ 0\}$ be a solution to  \eqref{eq:Schr2+}. 
Let $\chi \in C_0^{\infty}({\mathbb R})$, $\chi(t) = 1$ for $|t|\leq 1$, $\chi(t) = 0$ for $|t| \geq 2$ and define
$\chi_R(x) = \chi(|x|/R)$ for $R\geq 1$, $x \in {\mathbb R}^2_{+}$.

Suppose first that $\psi \in L^2({\mathbb R}^2_{+})$. 
One sees that $\nu \cdot (-i\nabla +{\bf \tilde F}) (\chi_R \psi) \big|_{\partial {\mathbb R}_{+}^2 }=0$, so
using \eqref{eq:Schr2+} and integration by parts
\begin{align}
\label{eq:IMS-Schr2+}
\langle \chi_R \psi , {\mathcal H} (\chi_R \psi) \rangle
= \lambda \| \chi_R \psi \|_2^2 + \frac{1}{R^2} \int_{{\mathbb R}^2_{+}} \big| \nabla \chi (x/R) \big|^2 \, |\psi (x)|^2\,dx.
\end{align}
Since $\psi \in L^2({\mathbb R}^2_{+})$ the last term in \eqref{eq:IMS-Schr2+} vanishes when $R \rightarrow \infty$. Therefore, using $\lambda < \Theta_0$ and the variational principle we obtain a contradiction to \eqref{eq:InfSpec}. So we conclude that
\begin{align}
\label{eq:NotL2}
\psi \notin L^2({\mathbb R}^2_{+}).
\end{align}

Clearly, \eqref{eq:NotL2} implies that $\| \chi_R \psi \|_{L^2({\mathbb R}^2_{+})} \rightarrow \infty$ as $R \rightarrow \infty$.
Notice that by the compact support of $\chi$ and the boundedness of $\psi$ we have, for some $C>0$ and all $R\geq 1$,
$$
\frac{1}{R^2} \int_{{\mathbb R}^2_{+}} \big| \nabla \chi (x/R) \big|^2 \, |\psi (x)|^2\,dx \leq C.
$$
(Here we used the fact that we study the $2$-dimensional problem.)\\
Thus, since $\lambda <\Theta_0$ by assumption,
\begin{align}
\frac{\langle \chi_R \psi , {\mathcal H} (\chi_R \psi) \rangle}{ \| \chi_R \psi \|_2^2}
= \lambda + o(1) < \Theta_0,
\end{align}
for $R$ sufficiently large. This is in contradiction to \eqref{eq:InfSpec} and thus $\psi$ cannot exist.
This finishes the proof for the equation \eqref{eq:Schr2+}.

We now prove the non-existence of bounded solutions to \eqref{eq:NL2+}. Let $\psi \in L^{\infty}({\mathbb R}^2_{+})\setminus \{ 0\}$ be a solution to  \eqref{eq:NL2+} and let the rest of the notation be as in the previous case. If $\lambda =0$ or $S=0$ the equation \eqref{eq:NL2+} is the same as  \eqref{eq:Schr2+}, so we may assume that $0<\lambda \leq \Theta_0$ and $S>0$.
Furthermore, after replacing $\psi$ by $S \psi$ we may assume that $S=1$.

Integrating by parts we obtain in analogy to \eqref{eq:IMS-Schr2+}
\begin{align}
\label{eq:IMS-NL2+}
\langle \chi_R \psi , {\mathcal H} (\chi_R \psi) \rangle
&= \lambda \| \chi_R \psi \|_2^2 
-\lambda \int_{{\mathbb R}^2_{+}} \chi_R^2 |\psi|^4\,dx\nonumber\\
&\quad\quad+ \frac{1}{R^2} \int_{{\mathbb R}^2_{+}} | \nabla \chi (x/R) |^2 \, |\psi (x)|^2\,dx.
\end{align}

Using \eqref{eq:InfSpec} in \eqref{eq:IMS-NL2+} we find
$$
\int_{{\mathbb R}^2_{+}} \chi_R^2 |\psi|^4\,dx \leq \lambda^{-1} \frac{1}{R^2} \int_{{\mathbb R}^2_{+}} | \nabla \chi (x/R) |^2 \, |\psi (x)|^2\,dx
\leq C,
$$
uniformly in $R$. Thus $\psi \in L^4({\mathbb R}^2_{+})$.
This we can use, as follows, to get a  good bound on the last term in \eqref{eq:IMS-NL2+},
\begin{align*}
\frac{1}{R^2} \int_{{\mathbb R}^2_{+}} | \nabla \chi (x/R) |^2 \, |\psi (x)|^2\,dx
\leq 
\frac{1}{R^3} \int_{{\mathbb R}^2_{+}} | \nabla \chi (x/R) |^4 \,dx
+ \frac{1}{R} \int_{{\mathbb R}^2_{+}} |\psi (x)|^4\,dx = o(1),
\end{align*}
as $R \rightarrow \infty$. Thus \eqref{eq:IMS-NL2+} becomes
\begin{align*}
\langle \chi_R \psi , {\mathcal H} (\chi_R \psi) \rangle
= \lambda \| \chi_R \psi \|_2^2  - \lambda \| \psi \|_4^4 + o(1)
< \lambda \| \chi_R \psi \|_2^2,
\end{align*}
for $R$ sufficiently large. Again this is a contradiction to \eqref{eq:InfSpec} and therefore implies that $\psi =0$.
\end{proof}

\subsection{Asymptotic estimates}~\\
In this and the following subsections we will use the non-existence results from Subsection~\ref{nonexist} to obtain improved versions of the estimates in Theorem~\ref{thm:Elliptic} in a reduced parameter range.
The application of this idea (`blow-up') to the Ginzburg-Landau system appeared to our knowledge first in \cite{LuPa1,LuPa2} and has since been used extensively since (see for instance \cite{LuPa5,PanCor,He-Pan}).

We consider solutions $(\psi, {\bf A})$ to 
\eqref{eq:GL} and satisfying the gauge condition \eqref{eq:gauge}.
Recall (see \eqref{infini}) that any solution of \eqref{eq:GL} satisfies the estimate $\| \psi \|_{\infty} \leq 1$.
If $\kappa/H$ is not too large, we can improve that estimate.

\begin{prop}\label{prop:1}~\\
Let $g: {\mathbb R}_+ \rightarrow {\mathbb R}_+$ satisfy that $g(\kappa) \rightarrow 0$ as $\kappa \rightarrow \infty$.
Then there exists a function $\tilde g$ with $\tilde g(\kappa) \rightarrow 0$ as $\kappa \rightarrow \infty$, such that if
$$
\kappa(\Theta_0^{-1} - g(\kappa)) \leq H  \leq \kappa(\Theta_0^{-1} + g(\kappa)),
$$
then any solution $(\psi, {\bf A})_{\kappa,H}$ of \eqref{eq:GL} satisfies
$$
\| \psi \|_{\infty} \leq \tilde g(\kappa).
$$
\end{prop}
\begin{remark}~\\
The upper bound $H \leq \frac{\kappa}{\Theta_0} + o(\kappa)$ is natural---at least in the study of minimizers of the Ginzburg-Landau functional ${\mathcal E}_{\kappa,H}$---for the following reason. It is known that for a given $\kappa$ there exists $H_{C_3}(\kappa)$ such that for all $H> H_{C_3}(\kappa)$, the unique minimizer of ${\mathcal E}_{\kappa, H}$ is the configuration $(\psi, {\bf A}) = (0, {\bf F})$ (up to change of gauge).
This {\it critical field} $H_{C_3}(\kappa)$ has been intensively studied \cite{Giorgi-Phillips,LuPa1, PiFeSt,He-Pan,FournaisHelffer3} and it is known to satisfy
$$
H_{C_3}(\kappa) = \frac{\kappa}{\Theta_0} + o(\kappa),
$$
for large $\kappa$ (actually much more precise asymptotic expansions exist).
With this notation it would be more natural, in the case of minimizers, to write the upper bound on $H/\kappa$ as
$H \leq H_{C_3}(\kappa)$.\\
Furthermore, it is known (see for instance \cite{LuPa1,FournaisHelffer}) that $\|\psi\|_{\infty}$ does not tend to zero if $H/\kappa \rightarrow \theta < \Theta_0$. Therefore the parameter domain of Proposition~\ref{prop:1} is optimal.
\end{remark}

The same type of argument as will be given in the proof of Proposition~\ref{prop:1} is used to prove the estimates below. These are slightly improved versions of \cite[Prop.~4.2] {He-Pan} and \cite[Lemma~7.1]{Panhc2}.

\begin{prop}
\label{prop:blowup}~\\
Let $0< \lambda_{\rm min} \leq \lambda_{\rm max}$.
There exist constants $C_0, C_1$ such that, if 
$$
\kappa \geq C_0,\quad\quad\quad\quad
\lambda_{\rm min} \leq \kappa/H \leq \lambda_{\rm max},
$$
then any solution $(\psi, {\bf A})$ of \eqref{eq:GL} satisfies
\begin{align}
\label{eq:first}
&\| p_{\kappa H {\bf A}} \psi \|_{C(\overline{\Omega})} \leq C_1 \sqrt{\kappa H} \| \psi \|_{\infty},\\
\label{eq:second}
&\| \curl {\bf A} -1 \|_{C^1(\overline{\Omega})}\leq C_1 \frac{1}{\sqrt{\kappa H}} \| \psi \|_{\infty}^2,\\
\label{eq:third}
&\| \curl {\bf A} -1 \|_{C^2(\overline{\Omega})}\leq C_1  \| \psi \|_{\infty}^2.
\end{align}
\end{prop}

\subsection{Extraction of convergent subsequences}~\\
The technique of proof of the estimates in Propositions~\ref{prop:1} and \ref{prop:blowup} is to study certain limiting equations. We will discuss this procedure here.

Let $\{P_n \}_{n} \subset \overline{\Omega}$ be a sequence of points and let 
$(\psi_n, {\bf A}_n)_{\kappa_n,H_n}$ be a sequence of solutions to \eqref{eq:GL} with $\psi_n \neq 0$.
We assume that $0<\lambda_{\rm min}  :=\liminf \kappa_n/H_n$,
$\limsup \kappa_n/H_n=:\lambda_{\rm max} < +\infty$.
We will proceed by repeatedly extracting subsequences of this original sequence. For convenience of notation we will not change the notation after each such extraction.
As detailed below, the result of this procedure will be that there exists a subsequence of the original sequence which (after rescaling and eventually a change of coordinates) converges to the solution of a limiting problem.

By extracting a subsequence (still indexed by $n$) we may assume that $P_n \rightarrow P \in \overline{\Omega}$. Similarly, we may assume that $\kappa_n/H_n \rightarrow \Lambda \in [\lambda_{\rm min},\lambda_{\rm max}]$. Also denote $S_n := \| \psi_n \|_{\infty} \neq 0$. We may assume that $S_n \rightarrow S \in [0,1]$.

By \eqref{eq:A2}, $\{{\bf A}_n\}_{n}$ is bounded in $W^{2,p}(\Omega)$, for all $p<\infty$. By compactness of the inclusion $W^{2,p}(\Omega) \hookrightarrow W^{s,p}(\Omega)$ for $s < 2$, we may extract a convergent subsequence (still denoted by ${\bf A}_n$).
Furthermore, for a given $\alpha <1$, we may choose $p$ sufficiently big and $s$ sufficiently close to $2$ in order to have the inclusion
$W^{s,p}(\Omega) \hookrightarrow C^{1,\alpha}(\overline{\Omega})$.
Thus we get the existence of some $\overline{\bf A} \in C^{1,\alpha}(\overline{\Omega})\cap W^{s,p}(\Omega)$ such that 
$$
{\bf A}_n \rightarrow \overline{\bf A} \quad \text{ in } \quad C^{1,\alpha}(\overline{\Omega}) \cap W^{s,p}(\Omega).
$$
We now identify the field generated by $\overline{\bf A}$.
The inequality \eqref{ineqimproved} holds for ${\bf A}_n$:
$$
\|  \curl {\bf A}_n-1 \|_2\leq C \frac{1}{ H_n}
 \|\psi_n\|_\infty\; \|\psi_n \|_2\;,
$$
with a constant $C$ independent of $n$ (only depending on $\Omega$). By passing to the limit (using \eqref{infini}), we find that
\begin{align}
\label{eq:overlineA}
\curl \overline{\bf A} = 1.
\end{align}

By passing to a subsequence we may assume that we are in one of the two cases below.

\noindent{\bf Case 1.} 
$$
\sqrt{\kappa_n H_n} \dist(P_n, \partial \Omega) \rightarrow \infty. 
$$
\noindent{\bf Case 2.} 
$$
\dist(P_n, \partial \Omega) \leq C/\sqrt{\kappa_n H_n},
$$
for some $C>0$.

\noindent{\bf Limiting equation for Case 1.}\\
Define, for any $R>0$ the following functions on the disc $B(0,R)$:
\begin{align*}
{\bf a}_n(y) &:= \frac{{\bf A}_n(P_n +y/\sqrt{\kappa_n H_n}) - {\bf A}_n(P_n)}{1/\sqrt{\kappa_n H_n}},\\
\varphi_n(y) &:= S_n^{-1} e^{-i\sqrt{\kappa_n H_n}{\bf A}_n(P_n)\cdot y} \psi_n(P_n +y/\sqrt{\kappa_n H_n}). 
\end{align*}
Since we are in Case 1, ${\bf a}_n$, $\varphi_n$ are defined on $B(0,R)$ for all $n$ sufficiently large.

Define the linear function $\tilde{\bf F}(y) := \big(D \overline{\bf A}(P)\big) y$.
By the convergencies 
$P_n \rightarrow P$,
${\bf A}_n \rightarrow \overline{\bf A}$ in $C^{1,\alpha}(\overline{\Omega})$, we find 
that
$$
{\bf a}_n \rightarrow \tilde{\bf F},
$$
in $C^{\alpha}\big(\overline{B(0,R)}\big)$ for all $R$.
By \eqref{eq:overlineA} we obtain
$$
\curl \tilde{\bf F} = 1, \quad\quad\text{ in } \quad\quad {\mathbb R}^2.
$$

The equation for $\psi$ in \eqref{eq:equationA} implies, since $\Div {\bf a}_n = 0$, that
\begin{align}
\label{eq:phi}
-\Delta \varphi_n - 2i {\bf a}_n \cdot \nabla \varphi_n + |{\bf a}_n|^2 \varphi_n
= \frac{\kappa_n}{H_n} (1 - S_n^2 |\varphi_n|^2)\varphi_n.
\end{align}
Notice that \eqref{eq:AF1} implies that for all $\alpha <1$, $\| {\bf a}_n \|_{C^{\alpha}(B(0,R))} \leq C_{\alpha}(R)$ for some $C_{\alpha}(R)>0$.
Also we have $\| \varphi_n \|_{\infty} \leq 1$.
Elliptic regularity (see \cite[Theorem~8.32]{GilbargTrudinger}) now implies, since $\frac{\kappa_n}{H_n}, S_n$ are bounded uniformly in $n$, the existence of a constant $C'_{\alpha}(R)>0$ such that
$$
\| \varphi_n \|_{C^{2,\alpha}(\overline{B(0,R/2)})}\leq C'_{\alpha}(R).
$$
Since the inclusion $C^{2,\alpha}(\overline{B(0,R/2)}) \hookrightarrow C^{2,\alpha'}(\overline{B(0,R/2)})$ is compact for any $\alpha' < \alpha$, we may for any $\alpha <1$, $R\geq 1$ extract a subsequence---denoted by $\{ \varphi_n^R\}$---having a limit in the $C^{2,\alpha}(\overline{B(0,R/2)})$ topology.
A `diagonal sequence' argument now gives the existence of a subsequence $\{ \tilde{\varphi}_n \}$ of the original sequence $\{ \varphi_n \}$ and a $\varphi \in C^{2,\alpha}({\mathbb R}^2)$ such that 
$$
\lim_{n\rightarrow \infty} \| \tilde{\varphi}_n - \varphi\|_{C^{2,\alpha}(\overline{B(0,R)})} = 0,
$$
for all $R>0$.
Passing to the limit in \eqref{eq:phi} we obtain the equation for $\varphi$:
\begin{align}
\label{eq:LimNoBdry}
(-i\nabla +{\bf \tilde F})^2 \varphi = \Lambda (1 - S^2 |\varphi|^2)\varphi.
\end{align}
\noindent{\bf Limiting equation for Case 2.}\\
The idea in the second case is the same as before but things are complicated slightly by the presence of the boundary. We make a change of variables in order to find a model on the half-plane.

Since we are in Case 2, $P \in \partial \Omega$. Let $Q_n \in \partial \Omega$ be the unique (for $n$ sufficiently large) boundary point such that $|P_n - Q_n| = \dist(P_n, \partial \Omega)$.
Let ${\mathcal O}$ be a (sufficiently small) neighborhood of $P$, let $\gamma : [-s_0,s_0] \rightarrow \partial \Omega$ be a smooth parametrization of the boundary with $\gamma(0)=P$, $|\gamma'(s)|=1$, and let $\nu(s)$ be the inward normal vector to $\partial \Omega$ at the point $\gamma(s)$. We may assume that $\{\gamma'(s),\nu(s)\}$ is a positively oriented basis.
Define the coordinate change 
$$
\Phi \;:\; (-s_0,s_0) \times (0,t_0) \rightarrow \Omega\cap {\mathcal O},
$$
by $\Phi(s,t) = \gamma(s)+ t\nu(s)$.
For $s_0,t_0,{\mathcal O}$ sufficiently small the map $\Phi$ is a diffeomorphism.

Let $\gamma_n$ be as $\gamma$ above, but with $\gamma_n(0) = Q_n$.
We now define $\Phi_n$ to be the same construction but with $\gamma$ replaces by $\gamma_n$ and $s_0$ replaced by $s_0/2$.
Since $Q_n \rightarrow P$ as $n \rightarrow \infty$ the image of $\Phi_n$ will contain $\Phi\big((-s_0/4,s_0/4) \times (0,t_0)\big)$ when $n$ is large.

Define 
\begin{align*}
&\tilde \psi_n := \psi_n \circ \Phi_n, &
&{\bf \tilde A}_n := (D \Phi_n)^t ({\bf A}_n \circ \Phi_n),\\
&J_n:= |\det D\Phi_n |,&
&M_n =\{M^n_{j,k}\}:= \big[ (D \Phi_n)^t (D \Phi_n) \big]^{-1}.
\end{align*}
Notice that $M_n\big|_{t=0} = Id$, and that the boundary condition $\nu\cdot {\bf A}_n \big |_{\partial \Omega}= 0$ implies that 
$$
e_2 \cdot {\bf \tilde A}_n\big |_{t=0} = 0.
$$

Implementing this change of variables in the equation \eqref{eq:equationA} for $\psi_n$ yields
\begin{align*}
&J_n^{-1} (-i\nabla + \kappa_n H_n {\bf \tilde A}_n)\cdot \big[ J_n M_n (-i\nabla + \kappa_n H_n {\bf \tilde A}_n) \tilde \psi_n \big]
= \kappa_n^2 (1 - |  \tilde \psi_n|^2) \tilde \psi_n,\\
&e_2 \cdot (-i\nabla + \kappa_n H_n {\bf \tilde A}_n) \tilde \psi_n \big|_{t=0} = 0.
\end{align*}
Let us calculate $\curl {\bf \tilde A}_n$. We use the geometric fact that $\nu_n'(s) = -k_n(s) \gamma_n'(s)$, where $k_n(s)$ is the curvature of the boundary at the point $\gamma_n(s)$.
Then 
$$
{\bf \tilde A}_n = (\tilde A_1^n, \tilde A_2^n) = 
\big((1-tk_n(s)) \gamma'_n(s) \cdot {\bf A}_n(\Phi_n(s,t)), \nu_n(s)\cdot {\bf A}_n(\Phi_n(s,t))\big).
$$
A direct calculation now yields
\begin{align}
\label{eq:curlAtilde}
\curl {\bf \tilde A}_n := \partial_s \tilde A_2^n - \partial_t \tilde A_1^n =
(1- tk(s)) (\curl {\bf A}_n)\big|_{\Phi_n(s,t)}.
\end{align}

Define $y_n :=\Phi_n^{-1}(P_n)$ and $z_n:=\sqrt{\kappa_n H_n} y_n$. Since we are in Case 2, $\{z_n\}$ is bounded and we may assume that $z_n \rightarrow z \in {\mathbb R}^2_{+}$.

We proceed to rescale as before. Define, with $\zeta=(\sigma,\tau)$,
\begin{align*}
{\bf a}_n(\zeta) &:= \frac{{\bf \tilde A}_n(\zeta/\sqrt{\kappa_n H_n}) - {\bf \tilde A}_n(0)}{1/\sqrt{\kappa_n H_n}},&
j_n(\zeta)&:=J_n(\zeta/\sqrt{\kappa_n H_n}),\\
\varphi_n(\zeta) &:= S_n^{-1} e^{-i\sqrt{\kappa_n H_n}{\bf \tilde A}_n(0)\cdot \zeta} \psi_n(\zeta/\sqrt{\kappa_n H_n}),
&
m_n(\zeta)&:=M_n(\zeta/\sqrt{\kappa_n H_n}). 
\end{align*}
We denote the components of ${\bf a}_n, m_n$ in the natural way, i.e.
${\bf a}_n=(a_1^n, a_2^n)$, $m_n = \{ m_{j,k}^n\}_{j,k=1}^2$.
Remember also the relations
$$
m_n\big |_{\tau=0} = Id,\quad\quad\quad\quad
e_2 \cdot {\bf a}_n \big |_{\tau=0} = 0.
$$

We get the resulting equation for the scaled function $\varphi_n$
\begin{align}
\label{eq:Border}
&j_n^{-1} (-i\nabla +  {\bf a}_n)\cdot \big[ j_n m_n (-i\nabla +  {\bf a}_n) \varphi_n \big]
= \frac{\kappa_n}{H_n} (1 - S_n^2 |  \varphi_n|^2) \varphi_n,\\
&e_2 \cdot (-i\nabla +  {\bf a}_n) \varphi_n \big|_{t=0} = 0.\nonumber
\end{align}

By \eqref{eq:A2} $\{{\bf A}_n\}$ and therefore $\{{\bf \tilde A}_n\}$ are bounded in $W^{2,p}$, $\forall p<\infty$. Thus $\{{\bf a}_n\}$ is bounded in $W^{1,p}\big( B(0,R) \cap {\mathbb R}^2_{+}\big)$ for all $R>0$.
We will below use standard results on elliptic regularity to conclude that
\begin{align}
\label{eq:NeumannElliptic}
\{\varphi_n\}_n \text{ is bounded in } W^{2,p}\big( B(0,R) \cap {\mathbb R}^2_{+}\big)
\text{ for all } R>0.
\end{align}

To prove \eqref{eq:NeumannElliptic} we rewrite the equation for $\varphi_n$ as follows.
\begin{align}
\label{eq:Extendable}
-\Div ( m_n \nabla \varphi_n) + {\bf b}_n \cdot \nabla \varphi_n + c_n \varphi_n = f_n,
\end{align}
with
$$
f_n:= \frac{\kappa_n}{H_n} (1 - S_n^2 |  \varphi_n|^2) \varphi_n +
i \big( \sum_{j,k} m^n_{j,k} \, \partial_j a_k^n\big) \varphi_n,
$$
and with the standard Neumann boundary condition $e_2 \cdot \nabla \varphi_n \big|_{\tau=0} =0$.
Here $\{ f_n\}$ is uniformly bounded in $L^p\big( B(0,R) \cap {\mathbb R}^2_{+}\big)$ (for all $R>0$) since $\| \varphi_n \|_{\infty} \leq 1$, and the coefficients ${\bf b}_n, c_n$ are uniformly bounded in 
$W^{1,p}\big( B(0,R) \cap {\mathbb R}^2_{+}\big) \hookrightarrow L^{\infty}\big( B(0,R) \cap {\mathbb R}^2_{+}\big)$.

In order to remove the boundary condition we extend by reflection. We denote extended functions by a superscript tilde. These functions will be defined by the fact that they are extensions of the original functions and that they are even or odd under the symmetry $(\sigma,\tau) \mapsto (\sigma,-\tau)$.
Those symmetry properties are as follows
\begin{align*}
&\tilde{\varphi}_n, \tilde{m}_{1,1}^n, \tilde{m}_{2,2}^n, \tilde{b}^n_1, \tilde{c}_n, \tilde{f} \quad \text{ are even},\\
&
\tilde{m}_{1,2}^n, \tilde{m}_{2,1}^n, \tilde{b}^n_2 \quad \text{ are odd}.
\end{align*}
Since $m_n \big|_{\tau=0} = Id$, the matrix $\tilde{m}_n$ thus defined is continuous and $\tilde{\varphi}_n$ satisfies the extended version of \eqref{eq:Extendable} (with symbols having a superscript tilde).
Clearly the bounded properties of ${\bf b}_n, c_n$ imply that
$\tilde{\bf b}_n, \tilde c_n$ are bounded in $L^{\infty}(B(0,R))$ for all $R$.
We can now apply the `interior' estimates \cite[Theorem~9.11]{GilbargTrudinger} to this extended equation and conclude that
\begin{align}
\| \varphi_n \|_{W^{2,p}( B(0,R) \cap {\mathbb R}^2_{+})}
&\leq \| \tilde \varphi_n \|_{W^{2,p}( B(0,R)}\nonumber\\
&\leq C\big( \| \tilde \varphi_n \|_{L^{p}( B(0,2R)}+  \| \tilde f \|_{L^{p}( B(0,2R)}\big)\nonumber\\
&\leq C' \big( \| \varphi_n \|_{L^{p}( B(0,2R)\cap {\mathbb R}^2_{+}}+  \|  f \|_{L^{p}( B(0,2R)\cap {\mathbb R}^2_{+}}\big).
\end{align}
Using that $\| \varphi_n \|_{\infty} \leq 1$ in order to get a uniform bound to $\| \varphi_n \|_{L^{p}( B(0,2R)\cap {\mathbb R}^2_{+}}$, we have therefore proved \eqref{eq:NeumannElliptic}.

With \eqref{eq:NeumannElliptic} established, we can proceed essentially as in the Case 1.
Let $\alpha<1$ and let $s<2$, $p<\infty$ be chosen such that 
$W^{s,p}\big( B(0,R) \cap {\mathbb R}^2_{+}\big) \hookrightarrow C^{1,\alpha}\big(\overline{ B(0,R) \cap {\mathbb R}^2_{+}}\big)$.
A diagonal sequence argument, as for Case 1, gives the existence of $\varphi \in W^{s,p}_{\rm loc}({\mathbb R}^2_{+}) \cap C^{1,\alpha}({\mathbb R}^2_{+})$ such that (eventually after extraction of a subsequence)
$$
\lim_{n\rightarrow \infty}
\| \varphi_n - \varphi \|_{C^{1,\alpha}\big(\overline{ B(0,R) \cap {\mathbb R}^2_{+}}\big)} = 0,
$$
for all $R>0$.
Furthermore, since $\| \varphi_n \|_{\infty} \leq 1$ for all $n$, the same inequality is true for $\varphi$.

Passing to the limit in \eqref{eq:Border} we obtain the equation for $\varphi$,
\begin{align}
\label{eq:LimitBorder}
&(-i\nabla +  {\bf \tilde F})^2 \varphi
= \Lambda(1 - S^2|  \varphi|^2) \varphi,\\
&e_2 \cdot (-i\nabla + {\bf \tilde F}) \varphi \big|_{t=0} = 0.\nonumber
\end{align}
Here we used \eqref{eq:curlAtilde} to conclude that the limiting vector field---which we denote by ${\bf \tilde F}$---satisfies $\curl {\bf \tilde F} = 1$, so the notation is consistent.

\subsection{Proofs of Proposition~\ref{prop:1} and Proposition~\ref{prop:blowup}}
\begin{proof}[Proof of Proposition~\ref{prop:1}]~\\
The proof goes by contradiction. If Proposition~\ref{prop:1} is false then there exist $\epsilon_0>0$ and a sequence $(\psi_n, {\bf A}_n)_{\kappa_n,H_n}$ of solutions to \eqref{eq:GL} such that
$(\Theta_0^{-1} - g(\kappa_n)) \leq H_n/\kappa_n  \leq (\Theta_0^{-1} + g(\kappa_n))$, $\kappa_n \rightarrow \infty$ and
$$
\| \psi_n \|_{\infty} \geq \epsilon_0.
$$
Choose $P_n \in \overline{\Omega}$ such that
$\| \psi_n \|_{\infty} = |\psi_n(P_n)|$.
We now proceed to extract subsequences as described above.
We may assume that either Case 1 or Case 2 is satisfied.
In Case 1, we find the limiting equation \eqref{eq:LimNoBdry} with $\Lambda = \Theta_0$ and $S\geq \epsilon_0$. 
Proposition~\ref{prop:NonExist} implies, since $\Theta_0 \leq 1$, that $\varphi \equiv 0$.
However, by assumption 
\begin{align}
\label{eq:NotZero}
|\varphi(0)| = \lim_{n\rightarrow \infty} |\tilde{\varphi}_n(0)| = \lim_{n\rightarrow \infty} \frac{| \psi_n(P_n)|}{\| \psi_n \|_{\infty}} = 1.
\end{align}
This is a contradiction, so we conclude that Case 1 cannot occur.

Since Case 1 cannot occur we necessarily find that Case 2 occurs. 
Thus the limiting equation becomes \eqref{eq:LimitBorder} with $\Lambda = \Theta_0$, $S\geq \epsilon_0$. By Proposition~\ref{prop:NonExist}, $\varphi \equiv 0$, but
$$
|\varphi(z)| = \lim_{n\rightarrow \infty} |\varphi_n(z_n)| = \lim_{n\rightarrow \infty} \frac{| \psi_n(P_n)|}{\| \psi_n \|_{\infty}} = 1.
$$
Thus Case 2 is also impossible and we conclude that Proposition~\ref{prop:1} is satisfied.
\end{proof}

Actually, using that the parameter regime in \eqref{eq:NL2} is larger than for the half-plane case \eqref{eq:NL2+} we realize that the above proof actually also implies the following result of independent interest.
\begin{prop}\label{prop:PsiLilleInt}~\\
Let $\epsilon_0,\epsilon_1>0$ be such that $0<\Theta_0 - \epsilon_1<1 - \epsilon_0$.
Then there exist $\kappa_0,C>0$ such that if
$(\psi,{\bf A})_{\kappa,H}$ is a solution to \eqref{eq:GL} with $\psi \neq 0$,
$$
\kappa>\kappa_0, \quad\quad \Theta_0 - \epsilon_1 \leq \kappa/H \leq 1 - \epsilon_0,
$$
and $P \in \overline{\Omega}$ is such that $|\psi(P)| = \| \psi \|_{\infty}$, then
$$
\dist(P, \partial \Omega) \leq \frac{C}{\sqrt{\kappa H}}.
$$
\end{prop}

\begin{proof}[Proof of Proposition~\ref{prop:blowup}]~\\
{\bf Proof of \eqref{eq:first}.}~\\
Suppose \eqref{eq:first} is wrong. Then there exists a sequence $(\psi_n, {\bf A}_n)_{\kappa_n,H_n}$ of solutions to \eqref{eq:GL}, and a corresponding sequence of points $\{P_n\} \subset \Omega$ such that
$$
\frac{|p_{\kappa_n H_n {\bf A}_n} \psi_n(P_n)|}{\sqrt{\kappa_n H_n} \| \psi_n\|_{\infty}} \rightarrow \infty.
$$
After extracting subsequences as before we find (along the converging subsequence)
$$
\lim_{n\rightarrow \infty} \frac{|p_{\kappa_n H_n {\bf A}_n} \psi_n(P_n)|}{\sqrt{\kappa_n H_n} \| \psi_n \|_{\infty}}
= | (-i\nabla - {\bf \tilde F}) \varphi(z)| < \infty,
$$
where $z=0$ in Case 1 and $z = \lim_{n\rightarrow \infty} \sqrt{\kappa_n H_n}\Phi_n^{-1}(P_n)$ in Case 2.
This yields a contradiction, so we conclude that \eqref{eq:first} is correct.

{\bf Proof of \eqref{eq:second}.}~\\
This inequality is a consequence of \eqref{eq:first}.
Remember that 
$$
\curl^2 {\bf A} := (\partial_{x_2} \curl {\bf A}, -\partial_{x_1} \curl {\bf A}).
$$
Thus, by the Ginzburg-Landau equation \eqref{eq:equationA} and \eqref{eq:first}
\begin{align}
\label{eq:RightForDeriv}
\| \nabla \curl {\bf A} \|_{\infty} = \| \curl^2 {\bf A} \|_{\infty}
= \frac{1}{\kappa H} \| \Re\{ \overline{\psi} \, p_{\kappa H {\bf A}} \psi \} \|_{\infty}
\leq C \frac{1}{\sqrt{\kappa H}} \| \psi \|_{\infty}^2.
\end{align}
This is \eqref{eq:second} for the derivatives.

Furthermore, since $\curl {\bf A} - 1 = 0$ on $\partial \Omega$ and $\Omega$ is bounded we get
$$
\| \curl {\bf A} - 1 \|_{\infty} \leq C \frac{1}{\sqrt{\kappa H}} \| \psi \|_{\infty}^2,
$$ 
from \eqref{eq:RightForDeriv} by integration. This finishes the proof of \eqref{eq:second}.

{\bf Proof of \eqref{eq:third}.}~\\
The proof of this inequality follows the same idea as the proof of the pair of inequalities \eqref{eq:first}-\eqref{eq:second}. One needs to take one extra derivative. We refer the reader to \cite{He-Pan} for details.
\end{proof}

\bibliographystyle{amsalpha}

\end{document}